\documentclass{svjour3}

\smartqed  % flush right qed marks, e.g. at end of proof
\usepackage{graphicx}

\usepackage[fleqn]{amsmath}
\usepackage{alltt}
\usepackage{amssymb}

\newcommand{\F}{\mathbb F}

\newcommand{\Tr}{{\rm Tr}}
\newcommand{\ad}{{\rm ad}}
\newcommand{\Gal}{{\rm Gal}}
\newcommand{\Aut}{{\rm Aut}}
\newcommand{\Stab}{{\rm Stab}}
\newcommand{\GL}{{\rm GL}}

\begin{document}
%\begin{frontmatter}

\title{On the Equivalence of Quadratic APN Functions \thanks{Research supported by the Claude
Shannon Institute, Science
Foundation Ireland Grant 06/MI/006 and the Irish Research Council for Science, Engineering and Technology}}

\author{Carl Bracken%\thanks{Research supported by
%Irish Research Council for Science, Engineering and Technology
%Postdoctoral Fellowship}%\thanks{Research supported by the Claude Shannon Institute, Science Foundation Ireland Grant 06/MI/006}
\and
Eimear Byrne%\thanks{Research supported by the Claude
%Shannon Institute, Science
%Foundation Ireland Grant 06/MI/006}
\and
Gary McGuire%\thanks{Research supported by the Claude
%Shannon Institute, Science
%Foundation Ireland Grant 06/MI/006}
\and
Gabriele Nebe}

\institute{E. Byrne \at
              School of Mathematical Sciences, University College Dublin, Ireland\\ 
              \email{ebyrne@ucd.ie}           %  \\
              \and
           C. Bracken \at
              School of Mathematical Sciences, University College Dublin, Ireland\\
              \email{carlbracken@yahoo.com}
              \and
           G. McGuire \at
           School of Mathematical Sciences, University College Dublin, Ireland\\ 
           \email{gary.mcguire@ucd.ie}
           \and
           G. Nebe \at
           Lehrstuhl D f\"ur Mathematik, RWTH Aachen University, 52056 Aachen, Germany\\
           \email{Gabriele.Nebe@rwth-aachen.de}
}

\date{Received: date / Accepted: date}

\maketitle

\begin{abstract}
\noindent 
Establishing the CCZ-equivalence of a pair of APN functions is generally quite difficult.
In some cases, when seeking to show that a putative new infinite
family of APN functions is CCZ inequivalent
to an already known family, we rely on computer calculation for small values of $n$. 
In this paper we present a method to 
prove the inequivalence of quadratic APN functions 
with the Gold functions.
%NEW, check this, we do not solve the computational problem but 
%proof that all infinitely many functions are inequivalent 
Our main result is that a quadratic function is CCZ-equivalent to the APN Gold function $x^{2^r+1}$
if and only if it is EA-equivalent to that Gold function.
As an application of this result, we prove that
a trinomial family of APN functions
that exist on finite fields of order $2^n$ where $n\equiv2\textrm{ mod }4$
are CCZ inequivalent to the Gold functions. The proof relies on some knowledge of the 
automorphism group of a code associated with such a function.



\end{abstract}

%\end{frontmatter}

\section{Introduction}

Let $L$ be a finite field.
A function $f : L \longrightarrow L$ is said to be \emph{almost perfect nonlinear} (APN)
if the number of solutions in $L$ of the equation
\begin{equation}\label{apneq}
f(x+a)-f(x)=b
\end{equation}
is at most 2, for all $a,b\in L$, $a\not=0$.
If the number of solutions of (\ref{apneq}) in $L$ is at most $\delta$, we say $f$ is
differentially $\delta$-uniform.
Thus APN is the same as differentially 2-uniform.
A differentially 1-uniform function is also called a perfect nonlinear function,
or a planar function; however,
these do not exist in characteristic 2 because in that case
 if $x$ is a solution of (\ref{apneq}), so is $x+a$.
In this paper we only consider finite fields of characteristic 2.

The classical example of an APN function is $f(x)=x^3$,
which is APN (over any field) because (\ref{apneq}) 
is quadratic.
These were generalized to the Gold functions $f(x)=x^{2^r+1}$,
which are APN over $\mathbb{F}_{2^n}$ if $(n,r)=1$.
%has the form $x^{2^r}a +xa^{2^r} + a^{2^{r+1}}+b=0$.
 
APN functions were introduced in \cite{N} by Nyberg, who defined them as the mappings with
highest resistance to differential cryptanalysis. Since then
many papers have been written on APN functions, although not many
different families of such functions are known.
For some time the list of known (extended affine) inequivalent APN functions
comprised only monomial functions and was conjectured to be complete.
Since 2006 several new families of non-monomial APN functions
%that are not CCZ-equivalent to any monomial APN function 
have been discovered.
Two binomial families are presented in
Budaghyan-Carlet-Leander \cite{BCL}.
Two infinite families, one of which generalizes the binomial family,
were discovered in \cite{BBMM}.
One of these families consists of the trinomials in Equation \eqref{f} below,
which we will study in Sections \ref{APN1} and \ref{compare}.

%Since the recent discoveries of new familes of APN functions, 
%the search for new families has been quite active.
An important aspect of this problem, after establishing the APN property
%of a class of polynomials, 
is to check that the functions are really new,
i.e. that they 
are {\em inequivalent} to the known APN families.
The notions of equivalence most pervasive in the current literature
are {\em extended affine} (EA) and {\em Carlet-Charpin-Zinoviev} (CCZ)
equivalence \cite{CCZ}. EA-equivalence is finer than CCZ-equivalence and is
usually somewhat
easier to establish.

Two functions $f,g : L \longrightarrow L$
are called EA-equivalent
if there exist affine permutations $A_1,A_2$ and an affine map $A$ such that
$g=A_1\circ f \circ A_2 + A$.
The differential uniformity of a function
is an invariant of EA-equivalence.
However, a bijective function is not necessarily EA-equivalent to its inverse,
even though they have the same differential uniformity.

%A more general notion of equivalence has been suggested in \cite{CCZ},
%which is referred to as CCZ-equivalence.
Two functions are called CCZ-equivalent if the graph
of one can be obtained from the graph of the other by an affine permutation
of the product space.
Differential uniformity and resistance to linear and differential attacks are 
invariants of CCZ-equivalence, and unlike EA-equivalence, any
permutation \emph{is} always CCZ-equivalent to its inverse.

In the instance that a function $f:L \longrightarrow L$ is {\em quadratic},
%i.e. all exponents in the non-zero terms in $f$ have binary weight 2,
the map $f(x+y)+f(x)+f(y)$ is bilinear.
Therefore, the problem of testing the APN property of $f$ is 
reduced to obtaining an estimate on the size of the kernel of the linear 
map $f(x+a)+f(x)+f(a)$. For this reason, most of the known
non-monomial APN functions are in fact quadratic.

It turns out that in the case of quadratic functions the problem of establishing CCZ
equivalence can sometimes be reduced to checking EA-equivalence.
Yves Edel has asked recently in some conference presentations whether any two quadratic APN
functions are CCZ-equivalent if and only if they are EA-equivalent.
The main result of this paper is a partial answer:
we prove that a quadratic APN function is
CCZ-equivalent to a Gold function if and only if it is EA-equivalent to that Gold function.
Up to now, proofs that CCZ-equivalence implies EA-equivalence have been
by lengthy brute force computations for specific functions; see
the proof of Theorem 4 in \cite{BCL}, for example.
Our result is more general, holding for any quadratic function, and the proof is
by different methods.
Our methods will involve a study of the automorphism group
of a code determined by a quadratic function. For the Gold functions, this
group is known and has been determined by Berger \cite{TB}. 
We  combine our main result with a study of the automorphism group
to show that a new
family of APN functions
found by Bracken-Byrne-Markin-McGuire 
are CCZ inequivalent to any Gold function.
This family is a subclass of that given in \cite{BBMM}
and has the following description.
Let $k$ and $s$ be  odd coprime integers,
let $b,c\in \mathbb{F}_{2^{2k}}$ with $c \notin  \mathbb{F}_{2^k}$, 
and $b$ a primitive element of $ \mathbb{F}_{2^{2k}}$. 
The polynomials of the form
\begin{equation}\label{f}
f_s(x)= bx^{2^s+1}+(bx^{2^s+1})^{2^k}+cx^{2^{k}+1}
\end{equation}
are APN on $\mathbb{F}_{2^{2k}}$.
Previously, these polynomials were
demonstrated in \cite{BBMM} to be inequivalent in general to 
$x^{2^r+1}$ by using a computer to show 
the result for $k=3$ and $k=5$.

This paper is organized as follows.
In Section \ref{equivAPNcodes} we discuss some general background, including
the important connections between an APN function  and a certain 
associated code.
Section \ref{APN} discusses the particular case of quadratic APN functions,
and introduces the property we use in this paper.
In Section \ref{GOLD} we present our main result, which proves
that any quadratic APN function that is CCZ-equivalent to a Gold function must
be EA-equivalent to that Gold function.
Section  \ref{APN1} proves some results about automorphisms of family of APN functions
in Equation (\ref{f}), and Section \ref{compare} applies the results of the paper to that family.

\section{Equivalence of APN~functions and codes.}\label{equivAPNcodes}

Throughout the paper we fix a finite field $K := \F _{2^n}$ of
characteristic 2.
Let $\Tr$ denote the absolute
trace map from $K$ to $\mathbb{F}_2$. We write $\F_2^{2^n}= \F_2 ^K$,
implicitly fixing an ordering of $K$.
To a function $f:K\to K$ we associate a linear code 
$C_f \leq \F_2^{2^n}= \F_2 ^K$ as 
$$C_f:= \{ c^f_{\alpha, \beta,\epsilon } \mid \alpha , \beta \in K,
\epsilon \in \F_2 \} $$
where
$$c^f_{\alpha,\beta ,\epsilon }: K\to \F_2,  x \mapsto
\Tr (\alpha x) + \Tr (\beta f(x)) + \epsilon .$$
It was first observed in \cite{CCZ} that the dual code of $C_f$ has minimum distance 6
if and only if $f$ is APN.
Also \cite[Thm. 6]{BBMM} (first stated by John Dillon in a talk given at Banff in 2006 and later in \cite{Dillon}) shows that two functions $f,g$ are
CCZ-equivalent if and only if the associated linear binary codes
$C_f$ and $C_g$  are equivalent.
Recall that two codes $C,D \leq \F_2^N$ are {\em equivalent},
if there is some permutation $\pi \in S_N$ of the
coordinate places with $\pi(C) = D$. Explicitly,
$$\pi(C) = \{ c^f_{\alpha, \beta,\epsilon }\pi  \mid \alpha , \beta \in K,
\epsilon \in \F_2 \}$$ 
where 
$$c^f_{\alpha,\beta,\epsilon }\pi = x \mapsto \Tr(\alpha x\pi)+\Tr(\beta f(x\pi))+\epsilon.$$ 
The {\em automorphism group} of a code $C$ is defined as
$$\Aut (C) := \{ \pi \in S_N \mid \pi(C) = C \} .$$

\begin{remark} \label{autos}
Identifying the places of the codes with the elements of $K$, we
obtain certain canonical permutation groups:
\begin{itemize}
\item[(a)] $E:=(K,+)$, the additive group of $K$, isomorphic to
$\mathbb{Z}_2^n$, acting regularly on $K$ by $a: K\to K, x\mapsto a+x $.
\item[(b)] $M:=K^*$, the multiplicative group of $K$, isomorphic to
$\mathbb{Z}_{2^n-1}$ acting on $K$
by $a: K \to K, x\mapsto ax $.
Note that $K^*$ fixes $0$ and acts regularly on $K\setminus \{ 0 \} $.
\item[(c)] $\Gamma := \Gal (K/\F_2) = \langle \sigma \rangle  \cong \mathbb{Z}_n$,
the Galois group of $K$ acting on $K$ as the {\em Frobenius automorphism} 
$\sigma: K\to K, x\mapsto x^2$.
\end{itemize}
\end{remark}

This paper mainly treats the important class of 
{\em quadratic} APN functions $f:K\to K $.
Recall that the polynomial $f \in K[x]$ is quadratic
if for any non-zero $k \in K$, the function 
$f(x+k)+f(x)+f(k)$ is a linearized polyomial in $x$, or equivalently,
if it is $\F_2$-linear.
%Recall that the polynomial $f = \sum _{m=0}^{2^n-1} p_m x^m $ 
%is {\em quadratic} if and only if $p_m\neq 0 \Rightarrow 
%m = 2^a+2^b$ for some $a,b  $. 
The family of trinomials (\ref{f}) is quadratic, as is $x^{2^r+1}$.

The following proposition is well known -- it states that the additive group of the field
is contained in the automorphism group of a quadratic function.

\begin{proposition}
Let $f:K \longrightarrow K$ be quadratic. Then $(K,+) \leq \Aut (C_f) $.
\end{proposition}

\begin{proof}
%For each $k \in K$, define
%$\pi_k: c^f_{\alpha,\beta,\epsilon}(x) \mapsto c^f_{\alpha,\beta,\epsilon}(x+k).$
Since $f$ is quadratic, for each $k \in K$ we may write
$L(x+k):=f(x+k)+f(x)+f(k) = \sum_i{k_i x^{2^i}}$ for some $k_i \in K$.
Using this and that fact that $\Tr(a)=\Tr(a^2)$ for each $a \in K$ we obtain: 
\begin{eqnarray*}
  \pi_k(c^f_{\alpha,\beta,\epsilon}(x)) & = & c^f_{\alpha,\beta,\epsilon}(x+k)\\
  & = & \Tr(\alpha(x+k))+ \Tr(\beta f(x+k)) + \epsilon \\
  & = & \Tr(\alpha x) + \Tr(\beta(L(x+k)+f(x)+f(k))) + \Tr(\alpha k) + \epsilon \\
  & = & \Tr(\alpha x) + \Tr(\beta L(x+k)) + \Tr(\beta f(x)) + \Tr(\beta f(k))+ \Tr(\alpha k) + \epsilon\\
  & = & \Tr(\alpha x) + \Tr(\sum_i (\beta k_i)^{2^{-i}} x) + \Tr(\beta f(x)) + \Tr(\beta f(k))+ \Tr(\alpha k) + \epsilon\\  
  & = & c^f_{\alpha', \beta, \epsilon'}(x)
\end{eqnarray*}
where $\alpha' = \alpha + \sum_i (\beta k_i)^{2^{-i}}$ and $\epsilon' = \epsilon + \Tr(\beta f(k))+ \Tr(\alpha k).$ It follows that the map
$$\pi_k: C_f \longrightarrow C_f :  c^f_{\alpha,\beta,\epsilon}(x) \mapsto c^f_{\alpha,\beta,\epsilon}(x+k)$$
is an automorphism of $C_f$.
\end{proof}
\qed

We recall some basic definitions from group theory. Further background reading may be read in \cite{Robinson}

\begin{definition}
   Let $G$ be a group and let $H,N$ be subgroups of $G$, with $N$ normal.
   \begin{enumerate}
   \item
       The {\em normalizer} of $H$ in $G$, denoted $N_G(H)$,
       is the subgroup of $G$ comprising all $g \in G$
       such that $gHg^{-1} = H$. 
       %We say that the elements of $N_G(H)$ conjugate $H$ into itself.
   \item
       The {\em centralizer} of $H$ in $G$, for which we write $C_G(H)$
       is the subgroup of $G$ comprising all $g \in G$ such that $ghg^{-1}=h$
       for all $h \in H$. 
   \item  
       If $N \cap H$ is the identity then the group $NH$ is called the {\em semi-direct product} 
       of $N$ and $H$ and we write $N:H$. 
      
   \end{enumerate}
\end{definition}

%\begin{definition}
%Let ${\mathcal A}:= N_{S_{2^n}}(K,+) $ denote the normalizer of $(K,+)$
%in the symmetric group.
%${\mathcal A} $ consists of all permutations $\pi  \in S_{2^n}$
%that conjugate the set $(K,+)$ into itself.
%\end{definition}

For the remainder, we will write ${\mathcal A}:= N_{S_{2^n}}(K,+)$ to denote
the normalizer of $(K,+)$ in the symmetric group on $2^n$ elements.

This normalizer ${\mathcal A}$ plays the key role in establishing 
EA-equivalence via Theorem \ref{EAequiv} below. 

%The structure of ${\mathcal A}$  is well known and follows from
%general facts in elementary group theory.
%Recall that the group theory notation $A:B$ denotes a semidirect product.

\begin{proposition}
The normalizer ${\mathcal A}$ of $(K,+)$ in the symmetric group is the full affine group. 
That is,
$$
{\mathcal A} =  (K,+) : \GL_n (\F_2)  \cong
(\mathbb{Z}_2^{n}) : \GL _n (\F_2) .$$
\end{proposition}

\begin{proof}
%Recall that $(K,+)$ acts regularly on $K$, so for each $k\in K$ there is 
%a unique $\pi _k \in E$ such that $\pi _k(0) = k$.
%In fact this permutation $\pi _k$ is the map 
Since conjugation is a group automorphism, we obtain a group homomorphism from the normalizer into the
automorphism group
$$ \kappa : {\mathcal A} \to \Aut(K,+) \cong \GL_n(\F_2) , \pi \mapsto ( e\mapsto \pi e \pi^{-1}) .$$
Clearly, the kernel of $\kappa$ is 
the centralizer $C_{S_{2^n}}(K,+) $ of $(K,+)$ in $S_{2^n}$ and so 
${\mathcal A} / C_{S_{2^n}}(K,+) \cong \GL _n (\F_2)$.
Now $(K,+)$ acts regularly on itself via $\pi_k: x\mapsto x+k$.
We claim that $C_{S_{2^n}}(K,+)  = (K,+)$. 
It is clear that $(K,+) \subseteq C_{S_{2^n}}(K,+)$, since $(K,+)$ is abelian.
To see the converse inclusion let $\theta \in C_{S_{2^n}}(K,+)$.
Composing $\theta $ with the
inverse of the permutation $\pi _{\theta (0)} \in (K,+)$ we may assume that 
$\pi (0 ) = 0$. By assumption, $\theta = \pi_k \theta \pi_{-k}$ for all $k \in K$ and hence $\theta(x) = \theta(x-k) + k $ for all $x,k \in K$. In particular this gives
$\theta(k) = k$ for all $k$, so that $\theta$ is the identity.
%Then the regularity of the action of $(K,+)$ implies that
%$\theta $ is the identity:
%For $k\in K$ the permutation $\pi $ commutes with 
%$\pi_k\in E$ hence
%$$\pi (k) = \pi (\pi_k(0)) = \pi_k(\pi (0)) = \pi_k(0) = k $$
We deduce that $C_{S_{2^n}}(K,+) = (K,+)$.
%The automorphism group of $E$ is clearly the full linear group
%$\GL_n (\F_2) = \GL _{\F_2} (K)$.
%This group permutes the elements of $K$ and hence is a subgroup
%of $S_K=S_{2^n}$. 
The elements in $\GL_n(\F_2)$ stabilize $0\in K$,
hence $(K,+)$ meets $\GL_n(\F_2)$ at the identity and we conclude that 
the normalizer ${\mathcal A}$ is the semidirect product
as given in the proposition.
\end{proof}
\qed

\begin{remark}
The group ${\mathcal A}$ is also the automorphism group,
${\mathcal A} = \Aut (C_0)$,
of the first-order Reed-Muller code (see \cite[Ch. 13, Sec. 9]{MWS})
$$C_0 = \{ c_{\alpha , 0, \epsilon } \mid \alpha \in K, \epsilon \in \F_2 \} .$$
\end{remark}

The next `folklore' result is an important reformulation of EA-equivalence
in terms of codes. To our knowledge a proof has not appeared in 
literature. We include a proof here for completeness.
%except in the preprint \cite{CN}.

\begin{theorem}\label{EAequiv} %\cite{CN}
%(TODO : find a better reference)
${\mathcal A}$ 
acts on $\{ C_f \mid  f:K\to K \}$.
Functions $f$ and $g$ are
EA-equivalent functions if and only if the
codes $C_f$ and $C_g$  are in the
same ${\mathcal A}$-orbit.
\end{theorem}

In the proof it will be convenient to work with generator matrices. 
Let $N:=2^n$ denote the length of the code $C_f$.
By a generator matrix $G$
for $C_f$
we mean a matrix $G$ with row-space $C_f$.
Choosing an
$\F _2$-basis $(b_1,\ldots ,b_n)$ of $K$ we obtain
a generator matrix of the form
$G
= \left(\begin{array}{c} 
{\bf 1} \\ G_0' \\ G_f \end{array} \right)$
where ${\bf 1} \in \{ 1 \} ^{1 \times N} $
denotes the row consisting of $1$ only and
$G_0', G_f \in \F_2^{n \times N} $ are defined by indexing the
columns with the elements of $K$ as
$$
(G_0')_{i,x} := \Tr (b_{i} x), \ \
(G_f)_{i,x} :=
\Tr (b_{i} f(x)).
\ \ (\star ) $$
%(In fact if $f$ is an APN~function then the rows of $G$ are linearly independent \cite{CCZ}).
\\
Note that 
$G_0 := \left(\begin{array}{c}
{\bf 1} \\ G_0' \end{array} \right) \in \F_2^{(n+1)\times N}$
is a generator matrix for the first-order Reed-Muller code $C_0$.

The following Lemma is of independent interest.

\begin{lemma}\label{equal}
Let 
 $f, g: K\to K$. 
Then  $C_f = C_g$ if and only if 
$g=A_1 \circ f + A $ for some affine permutation $A_1$
and some affine map $A$.
\end{lemma}

\begin{proof}
Let
$G
= \left(\begin{array}{c} 
{\bf 1} \\ G_0' \\ G_f \end{array} \right)$ and 
$G'
= \left(\begin{array}{c} 
{\bf 1} \\ G_0' \\ G_g \end{array} \right)$ be generator matrices of 
$C_f$ resp. $C_g$ as above.
Then $C_f = C_g$ if and only if $G$ and $G'$ have the same rowspace, 
if and only if there are $B_1 \in \GL _n(\F_2) $, $B\in \F_{2}^{n\times n}$,
$t\in \F_2^{n\times 1}$ such that 
$$G_g = B_1 G_f + BG_0' + t {\bf 1} .$$
By $(\star )$ above, this means that for all $x\in K$ and all $1\leq i \leq n$ 
$$ 
\Tr (b_{i} g(x)) = \Tr (B_1 b_i f(x)) + \Tr (Bb_i x) + t_i $$
and hence
$g=A_1 \circ f + A $ with $A_1 =( B_1^{\ad },0)$ and $A=(B^{\ad } , t)$
where $a^{\ad }$ denotes the adjoint linear map of $a$ with respect
to the trace bilinear form.
A similar calculation shows the converse.
\end{proof}
\qed

Proof of Theorem \ref{EAequiv}:
{\bf (1)} 
We first show that 
$$\{ C_f \mid  f:K\to K \} = \{ C\leq \F_2^K \mid C_0 \subseteq C,\dim (C) \leq 
2n+1  \},$$
and in particular that ${\mathcal A} = \Aut (C_0)$ acts on this set.
\\
The inclusion $\subseteq $ is clear. So let $C\leq \F_2^K$ be a code of dimension $\leq 2n+1$ that contains 
$C_0$ and let 
$$G = \left( \begin{array}{c} G_0 \\ G_1 \end{array}\right) \in \F_2^{(2n+1)\times N}$$ be a generator matrix of $C$.
Let $T \in \F_2^{n\times n}$ denote the Gram matrix of the basis
$(b_1,\ldots , b_n)$ of $K$  with respect to the trace bilinear form,
$$T_{i,j} = \Tr (b_i b_j ) .$$
Then $T\in \GL_n(\F_2)$ by the non degeneracy of the trace.
For $x\in K$ let $f_x$ denote the column of index $x$ of
$T^{-1} G_1$ and define $f:K\to K$ by 
$f(x):= \sum _{i=1}^d (f_x)_i b_i  \in K$ the
corresponding element in $K$.
Then $G_1 = G_f$ and hence $C = C_f$.
\\
%To this aim we write $A_1 = (a_1,t_1) $ for some $a_1\in \GL_n(\F_2) $ and
%$t_1 = A_1(0) \in  K$ and $A=(a,t)$ with $a\in \F_2^{n\times n}$ and $t\in K$.
%Then 
%$$A_1(f(x)) + A(x) = a_1(f(x)) + a(x) + t + t_1 .$$
%Replacing $t$ by $t+t_1$ we may assume that $t_1 = 0$ and hence
%$A_1 = (a_1,0)$ is linear. 
%If $\left(\begin{array}{c} G_0 \\ G_1 \end{array} \right) $ 
%is a generator matrix for $C_f$ as above  (so $G_1= G_f$)
%then a generator matrix of $C_{A_1f}$ is given by 
%$\left(\begin{array}{c} G_0 \\ a_1 G_1 \end{array} \right) $.
%Hence the codes are equal. 
%Similarly, replacing $f$ by $f+A$ with $(f+A)(x) = f(x) + a(x) + t $ 
%for some linear mapping $a:K\to K$ and $t\in K$, 
%yields equality between the codewords
%$c^{f+A}_{0,\beta,0} \in C_{f+A} $  and 
%$c^f_{\alpha,\beta,\epsilon} \in C_{f} $ where
%$$\alpha = a^{\ad} (\beta ), \ \mbox{ and } \epsilon = \Tr (\beta t) $$
%and $a^{\ad }$ denotes the adjoint linear map of $a$ with respect
%to the trace bilinear form.
%This is because 
%$$\Tr (\beta (f(x)+a(x)+t)) = 
%\Tr (\beta f(x) ) + \Tr (a^{\ad} (\beta ) x) + \Tr (\beta t ) .$$
%\\
{\bf (2)} 
Now let $f,g:K\to K$ be EA-equivalent, so there are 
affine permutations $A_1,A_2$ and an affine mapping  $A$ such that
$g = (A_1\circ f \circ A_2) + A$.
We have to show that $C_g$ and $C_f$ are in the same 
orbit under ${\mathcal A}$.
By Lemma \ref{equal} we may assume that $A_1 =1 $ and $A = 0$ and hence
that $g=f\circ A_2$  
for some $A_2 \in (K,+):\GL_{\F _2}(K) \cong {\mathcal A}$.
This means that $g(x) = f(A_2(x))$ and $A_2$ induces a permutation of
the places $x\in K$ that are in ${\mathcal A}$.
\\
{\bf (3)} Finally we prove the converse implication.
Assume that there is some $\pi \in {\mathcal A}$ such that
$C_f = \pi (C_g) $.
Let $G = \left( \begin{array}{c} {\bf 1} \\ G_0' \\ G_g \end{array}
\right) $ be
the generator matrix of
$C_g$ as above.
Then $\pi (C_g )=C_f$ has a generator matrix
$$G\pi = \left( \begin{array}{c} {\bf 1} \\ G_0' \pi  \\ G_g  \pi \end{array} \right) $$
obtained by
multiplying $G$ with the permutation matrix $\pi$ from the right.
Since $\pi $ fixes the code $C_0$, there is $A\in \GL_n(\F_2)$ and
$t\in \F_2^{n\times 1} $ such that
$$G_0' \pi = A G_0' + t {\bf 1} .$$
Therefore there are matrices
$t_1\in \F_2^{n\times 1},B_1\in \F_2^{n\times n},A_1 \in \GL_n(\F_2)$ s.t.
$$\left(\begin{array}{c|c|c} 
1 & 0 & 0 \\ \hline 
t & A & 0 \\ \hline
t_1 & B_1 & A_1  \end{array} \right) 
\left( \begin{array}{c} {\bf 1} \\ \hline G_0' \pi \\ \hline G_g \pi \end{array}
\right) = 
\left( \begin{array}{c} {\bf 1} \\ \hline G_0'  \\ \hline G_f  \end{array}
\right) $$
reading as
$$G_f = A_1 G_g \pi + B_1 G_0 ' + t_1 {\bf 1 } 
= G_{ A_1 \circ g \circ \pi + (B_1,t_1) } .$$
Since $\pi $ is an affine permutation, and $(B_1,t_1)$ is an affine mapping
this yields that
$f = A_1 \circ g \circ \pi + (B_1,t_1) $
is EA-equivalent to $g$.
\qed

\section{Quadratic APN functions}\label{APN}

We now consider quadratic APN functions $h$ satisfying the property that all regular elementary abelian subgroups of Aut$(C_h)$ are conjugate to $(K,+)$. 
We will show that such functions satisfy Edel's conjecture, 
i.e.,  that CCZ-equivalence for this family implies EA-equivalence.
In fact the APN property is not required in what follows. However, our interest in CCZ-equivalence 
is usually restricted to the class of APN functions.

\begin{theorem}\label{main}
Let $h$ be a quadratic function such that $E:=(K,+) \leq 
\Aut (C_h) =: H \leq S_{2^n} $. 
Assume that for all $\pi \in S_{2^n}$ 
$$\pi E \pi^{-1} \leq H \Rightarrow \mbox{ there is some }
h_{\pi }\in H \mbox{ such that } 
\pi E \pi ^{-1} = h_{\pi } E h_{\pi }^{-1} .$$
If a quadratic function $f$ is CCZ-equivalent to $h$ then
it is also EA-equivalent to $h$.
\end{theorem}

\begin{proof}
Since $f$ and $h$ are CCZ-equivalent, 
there is $\pi \in S_{2^n}$ such that $\pi (C_f) = C_h$.
The subgroup $E \leq \Aut(C_f) $ 
is hence conjugated to $\pi E \pi^{-1} \leq \Aut(C_h)$. 
By assumption this implies that $h_{\pi} ^{-1} \pi $ normalizes $E$, and
hence $h_{\pi} ^{-1} \pi \in N_{S_{2^n}} (E) = {\mathcal A}$ and
$h_{\pi} ^{-1} \pi (C_f) = h_{\pi }^{-1} (C_h) = C_h$.
By Theorem \ref{EAequiv} this means that the two functions are EA-equivalent.
\end{proof}
\qed

Since all regular elementary abelian subgroups are conjugate 
in $S_{2^n}$,
the following corollary is a reformulation of the theorem 
above and suggests one strategy to prove Edel's conjecture 
for arbitrary quadratic APN~functions.

\begin{corollary}\label{corconjug}
Let $h$ be a quadratic function such that all 
regular elementary abelian subgroups of $\Aut(C_h)$ are
conjugate to $(K,+)$.
Then all quadratic functions $f$ that are CCZ-equivalent to $h$ 
 are indeed EA-equivalent to $h$.
\end{corollary}

Thus Edel's conjecture for APN functions is proved under the stated hypothesis of 
Corollary \ref{corconjug}. 
We do not know any quadratic APN functions $h$ for which the above property
does not hold, i.e., for which $\Aut (C_h)$ contains more than one conjugacy class of regular elementary abelian subgroups. We checked that it holds for all known APN functions of degree up to 7.
Note that this is not true for arbitrary functions, for example, linear functions $f$ have
$\Aut(C_f)$ equal to the affine linear group, which usually has several different conjugacy classes of elementary abelian subgroups of order $2^n$.

%We checked that  the following functions $h$ satisfy the hypothesis in Corollary \ref{corconjug}.
%Note that we omitted the Gold functions; they will be treated in
%a more general context below.

%\begin{cor}
%Let $f$ be a quadratic APN~function that is CCZ-equivalent to one of the following APN functions $h$. Then $f$ and $h$ are EA-equivalent.
%\end{cor}

%Here we should include a list of known quadratic APN functions
%for which we checked this condition and which hence
%satisfy Edel's conjecture. \marginpar{ADD}

%All of degree 6:   All of degree 7: 
%These are too many, probably we should refer to the
%MAGMA database that is under construction.

\section{Quadratic functions equivalent to the Gold function}\label{GOLD}

Well understood examples of quadratic APN~functions are 
the Gold functions
$$g:K\to K, \ x \mapsto x^{2^r+1}$$
for a fixed positive integer $r$ satisfying $(r,n)=1$. 
The automorphism group ${\mathcal G}$ of $C_{g}$ contains some
obvious automorphisms: the additive group of the field,
the multiplicative group of the field, and the Galois automorphisms.
Results of Berger \cite{TB} show that this is the full automorphism group, i.e.,
$${\mathcal G}:=\Aut (C_{g}) \cong (K,+) : K^* : \Gal(K/\F_2) =
EM\Gamma $$
(in the notation of Remark \ref{autos}) 
of order 
$|{\mathcal G} | = |E| \cdot |M| \cdot |\Gamma | = 2^n(2^n-1) n $.
The proof uses the classification of finite simple groups.

%In particular the action of  ${\mathcal G}$ on $K$  is primitive and 
%$E=(K,+)$ is the unique minimal normal subgroup of this primitive
%permutation group (see \cite[7.2.6]{Robinson} for example). 
%It is therefore characteristic in ${\mathcal G}$, but indeed the 
%condition of Theorem \ref{main} is not quite the same that 
%$E$ being characteristic. 
%In fact one easily checks with Magma that the full affine group 
%$\mathbb{Z}_2^4:\GL_4(\F_2) $ has two conjugacy classes of regular elementary 
%abelian subgroups, though the normal subgroup $\mathbb{Z}_2^4$ is
%characteristic in this group.  

We recall some basic definitions.

\begin{definition}
    Let $G$ be a finite group and let $H$ be a subgroup of $G$. We say that $H$ is a $p$-subgroup
    of $G$ if $H$ has order $p^r$ for some positive integer $r$. $H$ is called a Sylow $p$-subgroup
    of $G$ if $r$ is the greatest positive integer such that $p^r$ divides $|G|$.
\end{definition}

The well-known second Sylow theorem states that all Sylow $p$-subgroups of a group $G$ are conjugate in $G$.
Since any subgroup that is normal in $G$ forms its own conjugacy class, as a direct consequence of this Sylow theorem 
we have that if $H$ is a normal Sylow $p$-subgroup of $G$ then it is the unique subgroup of $G$ of that order.

\begin{lemma}\label{lemkunique} 
$(K,+)$ is the unique subgroup of ${\mathcal G}$ that is isomorphic to
$\mathbb{Z}_2^n$.
\end{lemma} 

\begin{proof} 
This is clear, if $n$ is odd, since then $2^n$ is the largest
$2$-power in $|{\mathcal G}|$ and
$(K,+)$ is a Sylow 2-subgroup of ${\mathcal G}$, which must be unique since $(K,+)$ is normal in ${\cal G}$ and all such Sylow 2-subgroups are conjugate. \\
Assume now that $n=2k$ is even and
let $T \cong \mathbb{Z}_2^n$ be an 
elementary abelian subgroup  of ${\mathcal G}$. 
Then any $x\in T$ satisfies $x^2 =1$ so in particular $x^2\in (K,+)$. 
Therefore $T$ is 
%conjugated to a subgroup $T$ of 
a subgroup of $S:=(K,+):\langle \tau \rangle = \{ x\in {\mathcal G} \mid x^2\in (K,+) \} $,
where $\tau = \sigma ^{k}: z\mapsto z^{2^{k}} \in \Gal(K/\F_2) $ is the 
Galois automorphism of order 2. 
It is easy to check that the centralizer of $\tau $ in $S$ is isomorphic to 
$(\F_{2^k},+) \times \langle \tau \rangle$, which has order $2^{k+1}$.
Now consider the natural epimorphism 
$S \to S/(K,+) \cong \langle \tau \rangle$ and assume that the
elementary abelian subgroup $T\leq S$ is not contained in the kernel of 
this map (i.e. assume that $T$ is not equal to $K$).
Then there exists some $s \in (K,+)$ such that $s\tau \in T$.
Now $T$ is abelian and is generated by $s\tau$ and $T \cap K$. 
Therefore $T \cap K$ 
%is normal in $T$ of 
has index 2 in $T$ and so has order $|T|/2=2^{2k-1}$.
Let $\theta \in T \cap K$. Then $s \tau \theta = \theta s \tau = s \theta \tau$, and hence $\theta$ commutes with $\tau$. 
This shows that $T \cap K \subset C_S(\tau)$. But then
$2^{2k-1} = |T \cap K| \leq | C_S(\tau)|= 2^{k+1} <  2^{2k-1}$, giving a contradiction. 
We deduce that $T=(K,+)$, and hence $(K,+)$ is the unique elementary abelian
subgroup of order $2^{2k}$ of ${\mathcal G}$.
\end{proof}
\qed

The main result of this paper, stated below, follows now from 
Lemma \ref{lemkunique} and Corollary \ref{corconjug}.  

\begin{theorem}\label{maincor}
Let $f$ be a quadratic APN~function and $g$ be a Gold function. 
If $f$ and $g$ are CCZ-equivalent, then they are
EA-equivalent.
\end{theorem}

%Proof: combine the previous Lemma with Corollary \ref{corconjug}.

\begin{corollary}\label{edelhope}
Let $h$ be a quadratic APN~function such that 
 $\Aut(C_h) $ is isomorphic to a subgroup of $ {\mathcal G}$.
Then all 
quadratic APN~functions $f$ that are CCZ-equivalent to $h$ 
 are indeed EA-equivalent to $h$.
\end{corollary}

Regarding a proof of Edel's conjecture, we may indeed hope that the automorphism
group of any quadratic APN~function is contained in ${\mathcal G}$.
If this were true, proving it would complete a proof of Edel's conjecture, thanks
to Corollary \ref{edelhope}.
However, this is not true:

\begin{example}
Consider the quadratic functions given in \cite{Dillon} .
$$h_1:=x^3 + x^5 + u^{62}x^9 + u^3x^{10} + x^{18} + u^3x^{20} + u^3x^{34} + x^{40},$$
$$h_2:=x^3 + u^{11}x^5 + u^{13}x^9 + x^{17} + u^{11}x^{33} + x^{48} ,$$
and 
$$ h_3:=x^3 + x^{17} + u^{16}(x^{18} + x^{33}) + u^{15}x^{48}.$$
Then $h_1$ and $h_2$ are APN on $GF(2^6)$ and $|\Aut(C_{h_1})|=|\Aut(C_{h_2})| = 2^6.5$, which is not a divisor
of $2^6(2^6-1)6$. The polynomial $h_3$ is APN on $GF(2^8)$ and $\Aut(C_{h_3})$ has order $2^{10}.3^2.5$, which does not divide $2^8(2^8-1)8$.
\end{example} 

\section{Automorphisms of Family (\ref{f}) }\label{APN1}

We could now use Theorem \ref{maincor} directly to establish CCZ-inequivalence of a member of Family (\ref{f}) (or indeed any other quadratic) to the Gold functions by establishing EA-inequivalence, which can be achieved by a brute-force comparison of coefficients in the equation $g=A_1\circ f \circ A_2+A$. 

Instead we find that further knowledge of the automorphism group associated with Family (\ref{f}) allows us to show that for this family, CCZ-equivalence with the Gold functions holds not merely if and only if the corresponding codes are equivalent (EA-equivalence), but if and only if they are equal. Thus in this instance we can avoid applying brute-force. 

Let $k$, $s$ be odd coprime integers, 
$K = \F_{2^{2k}}$ and $L:=\F_{2^k}$ the subfield of $K$ of index 2. 
We denote by $T_2: K\to L$ the relative trace of $K$ to $L$.

We compute a subgroup ${\mathcal U}$ of the 
automorphism group of the APN functions in Family (\ref{f}), 
which is big enough to allow us to prove that if a function $f$ in Family  
(\ref{f}) is EA-equivalent to a Gold function $g$, then $C_f=C_g$.
%and hence there is a linear permutation $\varphi \in \GL_{n}(\F_2 )$,
%an endomorphism $\eta \in \End(K)$ and some $\tau \in K$ such that
%$$g(x) = \varphi (f(x) + \eta (x) + \tau )$$ for all $x\in K$.
We remark that the particular form of $f= T_2(bx^{2^s+1})+cx^{2^k+1}$ 
is helpful in determining some of the automorphisms of $C_f$.
Most other (known) APN functions do not have such a form, and determining
their automorphisms seems to be difficult.
 
It will be helpful to us to parametrize $f$ by $s$ and $c \in K \backslash L$; we write
$$f=f_{c,s}: = bx^{2^s+1}+(bx^{2^s+1})^{2^k}+cx^{2^{k}+1},$$
for any $b$ primitive in $K$. 

Since $f_{c,s}$ is an APN~function, 
$\dim (C_{f_{c,s}}) = 4k+1$ (c.f. \cite[Cor. 1]{CCZ}) and 
$$ C_{f_{c,s}} = \langle {\bf 1} \rangle \oplus C_0 \oplus C_c =
\langle c_{0,0,1} \rangle \oplus \{ c_{\alpha , 0, 0 } \mid \alpha \in K  \} \oplus \{ c^f_{0,\beta,0} \mid \beta \in K \}.$$

We claim the following.

\begin{lemma}\label{indep}
Any two $c,d \in K\setminus L$ define the same codes, i.e.,
$C_c = C_d$.
\end{lemma}

\begin{proof}  
For $c\in K\setminus L$ we have 
$$f_{c,s}(x) = (bx^{2^s+1}) + (bx^{2^s+1})^{2^k} + c x^{2^k+1} = 
T_2(bx^{2^s+1}) + c x^{2^k+1} .$$
Note that for any $x\in K$ the element $x^{2^k+1} $ lies in $L$.
We have to show that the set 
$\{ c^f_{0,\beta,0} \mid \beta \in K \}$ is independent of the 
choice of $c$.
By the transitivity of the trace we obtain 
$c^f_{0,\beta ,0}(x) =  \Tr _{L/\F_2} (T_2(\beta f_{c,s}(x))) $ and 
$$T_2( \beta f_{c,s}(x)) = T_2(\beta ) T_2(bx^{2^s+1}) + T_2(\beta c) x^{2^k+1} .$$
Since the trace $T_2$ is nondegenerate and 
$(1,c)$ as well as $(1,d)$ form a basis of $K$ over $L$, 
there is for any given pair
$(T_2(\beta ),T_2(\beta c))$ a unique $\beta ' \in K$ such that 
$$T_2(\beta ') = T_2 (\beta) \mbox{ and   } 
T_2(\beta ' d) = T_2 (\beta c)  .$$
So the code $C_{c}$ is independent of the choice of 
$c \in K\setminus L$.
\end{proof}
\qed

\begin{lemma} \label{F4star}  We have
$\mathbb{F}_4^* \subseteq \Aut(C_{f_{c,s}})$.
\end{lemma}

\begin{proof}
 Let $\omega$ be a generator for $\mathbb{F}_4^*$.
Since $s$ and $k$ are odd, the exponents $2^s+1$ and $2^k+1$ are 
both multiples of 3 and hence $f_{c,s}(\omega x)=f_{c,s}(x)$.
\end{proof}
\qed

\begin{lemma}
 We have $L^* \subseteq \Aut(C_{f_{c,s}})$.
\end{lemma}

\begin{proof}
If $z\in L^*$ it is easy to check that
$f_{c,s}(zx)=z^{2^s+1} f_{cz^{1-2^s}} (x)$.
All transformations involved do not change the code $C_{f_{c,s}}$,
using Lemma \ref{indep}.
So multiplication by a primitive element of $L$ is an automorphism.
\end{proof}
\qed

Hence we obtain the following result:

\begin{theorem}
$\Aut(C_{f_{c,s}})$ 
contains a subgroup ${\mathcal U} \cong (K,+) : (\mathbb{F}_4^* \times L^*) $ of order
$2^{2k} \cdot 3(2^k-1) $. 
\end{theorem}

Note that $\Aut(C_{f_{c,s}})$ is not abelian, since the subgroup ${\mathcal U}$
we know about is not abelian.

For $s=1$ we obtain one more automorphism giving rise to a subgroup
of order
$2^{2k} \cdot 3k \cdot (2^k-1)$ of $\Aut (C_{f_{c,1}})$.
We conjecture that this is the actual order. 
This has been verified by computer for $k=3$ and $k=5$.

\begin{lemma}
$\Aut(C_{f_{c,1}})$ has an element $\delta $ of order $3k$, such 
that $\delta ^k=\omega $ from Lemma \ref{F4star}.
\end{lemma}

\begin{proof}
Choose $c=b^{(2^k+1)/3}$.
Then $c\in K \setminus L$ and by Lemma \ref{indep} we may
assume without loss of generality that $f=f_{c,1}$. 
It is easily checked that $f_{c,1}(bx)$ is equal to
$\sigma^2 (f(x))$, where $\sigma$ is the Frobenius automorphism of $K$ 
over $\F_2$.
Letting $T_b(f(x))=f(bx)$, the map
$\delta := \sigma^{-2} \circ T_b$ is hence an automorphism of $C_{f_c}$,
and its order can be checked to be  $3k$.
\end{proof}
\qed

The automorphism groups of other families of APN functions do not appear to be
as easy to work with as for Family (\ref{f}).
Therefore we have  not been able to prove similar results for other families.

\section{Inequivalence.}\label{compare}

We apply the results of the previous sections to give a proof 
that Family (\ref{f}) functions are not CCZ-equivalent to Gold functions.

\begin{theorem}\label{ineq}
Let $g:K\to K$ be a Gold-function and $f:K\to K$ be an APN~function in
Family (\ref{f}). 
If $f$ and $g$ are EA-equivalent, then the associated codes 
$C_f$ and $C_g$ are equal.
\end{theorem}

For the proof of the theorem we need two lemmas, the first one is 
surely well-known. Recall that if a group $G$ acts on a set $X$ and $a \in X$
then the {\em stabilizer subgroup} of $a$, denoted $\Stab_G(a)$ is the set of
elements of $G$ that fix $a$. 

\begin{lemma}\label{norm}
$N_{\GL_n(\F_2)} (K^*) = K^* : \Gal (K/\F _2) $.
\end{lemma}

\begin{proof}
Let $G=N_{\GL_n(\F_2)} (K^*)$.
Clearly $K^* \leq G = K^* \Stab_G(1)$.
So it is enough to show that the elements $\pi \in G$ with
$\pi (1) = 1$ are indeed field automorphisms of $K$ and therefore
contained in $\Gal(K/\F_2)$. 
Choose $\pi \in G$ such that $\pi(1) =1$.
Since $\pi \in \GL_n(\F_2)$, the mapping $\pi $ acts linearly on 
the set $K$ and hence respects the addition.
We now show that $$\pi (ab) = \pi (a) \pi (b)
\mbox{ for all } a,b \in K .$$
To see this let $\alpha, \beta \in K^* \subset S_n$ be 
such that $\alpha (1) = a, \beta (1) = b $.
Since $\pi $ normalizes $K^*$, also
$$\tilde{\alpha } := \pi \alpha \pi ^{-1} \mbox{ and } 
\tilde{\beta } := \pi \beta \pi ^{-1}  \in K^* .$$
We calculate $\pi (a) = \pi (\alpha (1)) = (\pi \alpha \pi^{-1})(1) =
\tilde{\alpha }(1) $ and similarly $\pi(b) = \tilde{\beta }(1)$.
Clearly $ab = \alpha (\beta (1))$ and 
$$\pi (ab) = (\pi \alpha \beta \pi ^{-1} ) (1) =
\tilde{\alpha } (\tilde{\beta } (1)) =\pi (a) \pi (b) .$$
%So we have proven that $\pi $ acts as a field automorphism on $K$,
\end{proof}
\qed

\begin{lemma}\label{uniquecyc}
The group $G$ from Lemma \ref{norm} contains
%$G:= K^* : \Gal (K,\F_2) \cong \mathbb{Z}_{2^{2k}-1}:\mathbb{Z}_{2k} $
a unique cyclic subgroup of order $3(2^{k}-1)$.
\end{lemma}

\begin{proof}
All elements of $G$ are of the form $a \gamma $ where $a\in K^*$ and 
$\gamma \in \Gal (K,\F_2)$. 
Assume that such an element $a\gamma $ has order $3(2^k-1)$.
Let $\ell $ be the order of $\gamma $. Then $\ell $ divides $2k= \ell m $ and 
also the order of $a\gamma $. We calculate
$$(a\gamma )^{\ell } = N_{K/\F_{2^{m}}}(a) \in \F_{2^{m}}^*
 \cong \mathbb{Z}_{2^{m}-1}.$$
So $a\gamma $ has order dividing $\ell (2^{m} -1 )$ and we conclude that
$$3(2^k-1) \mbox{ divides } \ell  (2^{2k/\ell } - 1 ) $$ 
from which we obtain that 
$\ell  = 1$. 
%(TODO: really check this elementary calculation !! )
\end{proof}
\qed

\begin{proof}(of Theorem \ref{ineq})
Assume that there is some $\pi \in {\mathcal A}$ with 
$\pi (C_f) = C_g$. 
We identify the places of the code with $K$.
Since $\Aut(C_g)$ is 2-transitive on $K$ (cf. \cite{MWS}), 
we may assume without loss of generality that 
$$\pi (0) = 0 \mbox{ and } \pi (1) = 1 \mbox{ so } \pi \in \GL_n(\F_2) =\Stab _{{\mathcal A}} (0).$$
Moreover $\pi $ conjugates 
 ${\mathcal U} \leq \Aut (C_f)$ into ${\mathcal G} = \Aut (C_g)$, 
and since $\pi $ fixes $0$, also 
$$ \mathbb{Z}_{3( 2^{k}-1)} \cong \pi \Stab _{{\mathcal U}} (0) \pi ^{-1} 
= \pi (L^* \times \F_4^*) \pi^{-1} \leq 
\Stab _{{\mathcal G}} (0) =  K^* : \Gal (K,\F_2) .$$
By Lemma \ref{uniquecyc} this implies that $\pi $ normalizes 
$L^* \times \F_4^* \leq K^*$.
Since $L^*$ and $\F_4^*$ generate $K$ as an $\F_2$-algebra,
the linear span of the matrices in $L^*\times \F_4^*$ is equal to 
$K= K^* \cup \{ 0 \} \subset \F_2^{n\times n}$. 
Therefore $\pi $ also normalizes $K$ and hence $K^*$,
so $\pi \in \Gal (K,\F_2) \leq {\mathcal G}$ by Lemma \ref{norm}.
This proves the theorem since we have shown that any equivalence 
$\pi $ between the codes $C_f$ and $C_g$ is indeed already contained
in $\Aut (C_g)$.
\end{proof}
\qed
   
\begin{theorem}\label{unequal}
$C_f \neq C_g$.
\end{theorem}

\begin{proof}
Suppose that $C_f=C_g$. Then given any $\epsilon \in \F_2$, $\alpha,\beta \in K$ there exist $\epsilon' \in \F_2$, $\alpha',\beta' \in K$ satisfying
$$\epsilon +  \Tr(\alpha x)+\Tr(\beta(T_2(bx^{2^s+1})+ cx^{2^k+1})) = 
\epsilon' +  \Tr(\alpha' x)+\Tr(\beta'x^{2^r+1}))$$
for all $x \in K$, so in particular we must have $\epsilon =\epsilon'$. 
Choose $\beta \in L$. Then we have 
$\Tr(\beta(T_2(bx^{2^s+1}))) = \Tr(T_2(\beta b x^{2^s+1}))=0$ and so
$$ \Tr((\alpha +\alpha')x)=\Tr(\beta cx^{2^k+1}) +\Tr(\beta'x^{2^r+1}),$$
for all $x \in K$. Using the linearity of the LHS we obtain
$$\Tr(\beta c(x^{2^k}a+x a^{2^k} )+\beta'(x^{2^r}a+ x a^{2^r})) =
\Tr((\beta  x^{2^k}(c+c^{2^k}) + (\beta')^{2^{-r}}x^{2^{-r}}+\beta' x^{2^r})a)=0, $$
for all $x,a \in K$. This implies that $\beta  x(c+c^{2^k}) + (\beta')^{2^{k-r}}x^{2^{k-r}}+(\beta')^{2^{k}} x^{2^{k+r}} \in K[x]$ is identically zero, which, observing the degree of this polynomial and the fact that $(r,2k)=1$, we see is impossible unless $\beta = \beta'=0$. 
\end{proof}
\qed

\begin{remark} In fact this can also be readily seen by Lemma \ref{equal} 
by a simple comparison of coefficients.
\end{remark}

We now combine the results of Theorems \ref{ineq}, \ref{unequal} and Corollary \ref{maincor} in the following statement. 

\begin{corollary}
The functions of Family \ref{f} are not CCZ-equivalent to the Gold functions. 
\end{corollary}

\begin{acknowledgement} We thank the referees whose comments led to a much better presentation of this paper.
\end{acknowledgement}

\end{document}